\documentclass[12pt,reqno]{amsart}
\usepackage{amsmath, amsfonts, amssymb, amsthm}
\def\Z{\mathbb Z}

\def\1{{\bf 1}}

\def\jacob #1#2{\genfrac{(}{)}{}{}{#1}{#2}}

\def\qbinom #1#2#3{{\genfrac{[}{]}{0pt}{}{#1}{#2}}_{#3}}
\theoremstyle{plain}
\newtheorem{theorem}{Theorem}

\theoremstyle{definition}

\theoremstyle{remark}

\begin{document}
\title{The Rodriguez-Villegas type congruences for\\[5pt] truncated $q$-hypergeometric functions}
\author{Victor J. W. Guo}
\address{Department of Mathematics, Shanghai Key Laboratory of PMMP, East China Normal University,
500 Dongchuan Rd., Shanghai 200241, People's Republic of China}
\email{jwguo1977@aliyun.com}
\author{Hao Pan}
\address{Department of Mathematics, Nanjing University, Nanjing 210093, People's Republic of China}
\email{haopan79@zoho.com}
\author{Yong Zhang}
\address{Department of Basic Course, Nanjing Institute of Technology, Nanjing 211167,
People's Republic of China}
\email{yongzhang1982@163.com}
\maketitle

\section{Introduction}
\setcounter{lemma}{0}\setcounter{theorem}{0}\setcounter{corollary}{0}
\setcounter{equation}{0}

Define the truncated hypergeometric function
$$
{}_{2}F_1\bigg[\begin{matrix}x_1,\,x_2\\y_1\end{matrix}\bigg|z\bigg]_n
=\sum_{k=0}^{n-1}\frac{(x_1)_k(x_{2})_k}{(y_1)_{k}}\cdot\frac{z^k}{k!},
$$
where $(x)_k=x(x+1)\cdots(x+k-1)$ if $k\geq 1$ and $(x)_0=0$.
Motivated by the Calabi-Yau manifold, Rodriguez-Villegas \cite{RV03} conjectured some congruences on truncated hypergeometric functions modulo $p^2$ and $p^3$.
Nowadays, most of those conjectures have been confirmed. For example, with help of the Gross-Koblitz formula, Mortenson \cite{M03}
proved that, for any prime $p\geq 5$,
\begin{align}
{}_{2}F_1\bigg[\begin{matrix}1/2,\,1/2\\1\end{matrix}\bigg|1\bigg]_p &\equiv\jacob{-1}p\pmod{p^2},
\label{mor2} \\[5pt]
{}_{2}F_1\bigg[\begin{matrix}1/3,\,2/3\\1\end{matrix}\bigg|1\bigg]_p &\equiv\jacob{-3}p\pmod{p^2},
\label{mor3} \\[5pt]
{}_{2}F_1\bigg[\begin{matrix}1/4,\,3/4\\1\end{matrix}\bigg|1\bigg]_p &\equiv\jacob{-2}p\pmod{p^2},
\label{mor4} \\[5pt]
{}_{2}F_1\bigg[\begin{matrix}1/6,\,5/6\\1\end{matrix}\bigg|1\bigg]_p &\equiv\jacob{-1}p\pmod{p^2},
\label{mor6}
\end{align}
where $\jacob{\cdot}{p}$ is the Legendre symbol modulo $p$. Z.-W. Sun \cite{S13} gave an elementary proof for
\eqref{mor2}--\eqref{mor6}. Subsequently, Z.-H. Sun \cite{S14} generalized the above congruences to the following unified form:
\begin{equation}\label{szh}
{}_{2}F_1\bigg[\begin{matrix}\alpha,\,1-\alpha\\1\end{matrix}\bigg|1\bigg]_p\equiv(-1)^{\langle -\alpha\rangle_p}\pmod{p^2}.
\end{equation}
Here $\alpha$ is a rational number whose denominator is prime to $p$, and $\langle x\rangle_p$
denotes the integer in $\{0,1,\ldots,p-1\}$ such that $x\equiv \langle x\rangle_p\pmod{p}$. Note that we can also define
$\langle x\rangle_n$ similarly on condition that the denominator of $x$ is prime to $n$.

It is natural to define the truncated $q$-hypergeometric function as follows:
$$
{}_2\phi_1\bigg[\begin{matrix}x_1,\,x_2\\y\end{matrix}\bigg|q,z\bigg]_n=\sum_{k=0}^{n-1}\frac{(x_1;q)_k(x_2;q)_k}{(y;q)_{k}(q;q)_k}z^k,
$$
where
$$
(x;q)_k=\begin{cases}(1-x)(1-xq)\cdots(1-xq^{k-1}),&\text{if }k\geq 1,\\[5pt]
1,&\text{if }k=0.
\end{cases}
$$
Recently, Guo and Zeng \cite{GZ14} have obtained a $q$-analogue of \eqref{mor2}:
\begin{equation}\label{qmor2}
{}_{2}\phi_1\bigg[\begin{matrix}q,\,q\\q^2\end{matrix}\bigg|q^2,1\bigg]_p\equiv\jacob{-1}pq^{\frac{1-p^2}{4}}\pmod{[p]^2},\qquad p\geq 3,
\end{equation}
where $[p]=1+q+\cdots+q^{p-1}$ and the above congruence is considered over the polynomial ring $\Z[q]$.
Furthermore, they also conjectured that, for $p\geq5$,
\begin{align}
{}_{2}\phi_1\bigg[\begin{matrix}q,\,q^2\\q^3\end{matrix}\bigg|q^3,1\bigg]_p &\equiv\jacob{-3}pq^{\frac{1-p^2}{3}}\pmod{[p]^2},
\label{qmor3} \\[5pt]
{}_{2}\phi_1\bigg[\begin{matrix}q,\,q^3\\q^4\end{matrix}\bigg|q^4,1\bigg]_p &\equiv\jacob{-2}pq^{\frac{3(1-p^2)}{8}}\pmod{[p]^2},
\label{qmor4} \\[5pt]
{}_{2}\phi_1\bigg[\begin{matrix}q,\,q^5\\q^6\end{matrix}\bigg|q^6,1\bigg]_p &\equiv\jacob{-1}pq^{\frac{5(1-p^2)}{12}}\pmod{[p]^2}.
\label{qmor6}
\end{align}
In this paper, we shall prove the congruences \eqref{qmor3}--\eqref{qmor6}. More precisely, we shall give a $q$-analogue of (\ref{szh})
as follows.
\begin{theorem}\label{qrvt}
Let $n,d\geq 2$ with $\gcd(n,d)=1$ and let $r$ be an integer. Then
\begin{equation}\label{qrv}
{}_2\phi_1\bigg[\begin{matrix}q^r,\, q^{d-r}\\q^d\end{matrix}\bigg|q^d,1\bigg]_n\equiv(-1)^{a}q ^{(ad+r)(a-\frac{n-1}{2})-d\binom{a+1}2}\pmod{\Phi_n(q)^2},
\end{equation}
where $\Phi_n(q)$ denotes the $n$-th cyclotomic polynomial in $q$ and $a=\langle -r/d\rangle_n$.
\end{theorem}
For example, letting $r=1$, $d=3$ and letting $n=p$ be a prime greater than $3$, we have
$$a=\langle -1/3\rangle_p=\begin{cases}
(p-1)/3,&\text{if }p\equiv1\pmod{3},\\[5pt]
(2p-1)/3,&\text{if }p\equiv2\pmod{3},
\end{cases}$$
and so
$$
{}_2\phi_1\bigg[\begin{matrix}q,\,q^{2}\\ q^3\end{matrix}\bigg|q^3,1\bigg]_{p}\equiv\jacob{-3}{p}
q^{(ad+r)(a-\frac{n-1}{2})-d\binom{a+1}2}=\jacob{-3}{p}q^{\frac{1-p^2}{3}}\pmod{[p]_q^2},
$$
which is the congruence \eqref{qmor3}.

\section{Proof of Theorem \ref{qrvt}}
\setcounter{lemma}{0}\setcounter{theorem}{0}\setcounter{corollary}{0}
\setcounter{equation}{0}

Recall that the $q$-binomial coefficients ${n\brack k}_q$ are defined by
$$
\qbinom{n}{k}q=\frac{(q^{n-k+1};q)_k}{(q;q)_k},
$$
and the $q$-integer $[n]_q$ is defined as $[n]_q=\frac{1-q^n}{1-q}$.
Our proof of Theorem \ref{qrvt} only requires the following two forms of the $q$-Chu-Vandemonde identity (see, for example, \cite[(3.3.10)]{Andrews}):
\begin{align}
\qbinom{n+m}{k}q &=\sum_{j=0}^kq^{(n-j)(k-j)}\qbinom{n}{j}q\qbinom{m}{k-j}q, \label{eq:qChu-1} \\[5pt]
\qbinom{n+m}{k}q &=\sum_{j=0}^kq^{j(m-k+j)}\qbinom{n}{j}q\qbinom{m}{k-j}q.  \label{eq:qChu-2}
\end{align}

It is easy to see that
$$
\frac{(q^r;q^d)_k}{(q^d;q^d)_k}=(-1)^kq^{rk+d\binom{k}2}\qbinom{-r/d}{k}{q^d}.
$$
So writing $\alpha=-r/d$, the congruence \eqref{qrv} is equivalent to
\begin{align}
\sum_{k=0}^{n-1}q^{dk^2}\qbinom{\alpha}{k}{q^d}\qbinom{-1-\alpha}{k}{q^d}
\equiv (-1)^{a}q^{(\alpha-a)d(a-\frac{n-1}2)-d\binom{a+1}2}\pmod{\Phi_n(q)^2}.
\end{align}
Note that $a=\langle\alpha\rangle_n$. Let $s=(\alpha-a)/n$. Then $sd$ is an integer. By the $q$-Chu-Vandemonde identity \eqref{eq:qChu-1},
we have
\begin{align*}
\qbinom{a+sn}{k}{q^d}=&\sum_{j=0}^{k}q^{d(sn-j)(k-j)}\qbinom{sn}{j}{q^d}\qbinom{a}{k-j}{q^d}\\
\equiv&
\qbinom{a}{k}{q^d}-\sum_{j=1}^{k}(-1)^{j}q^{-dj(k-j)-d\binom{j}2}\frac{[sn]_{q^d}}{[j]_{q^d}}\qbinom{a}{k-j}{q^d}\pmod{\Phi_n(q)^2},
\end{align*}
where we have used the fact that
\begin{align*}
&\frac{[sn]_{q^d}}{[j]_{q^d}}=\frac{1-q^{sdn}}{1-q^{jd}}\equiv 0 \pmod{\Phi_n(q)}, \\[5pt]
&\qbinom{sn-1}{j-1}{q^d}\equiv\qbinom{-1}{j-1}{q^d}
=(-1)^{j-1}q^{-d\binom{j}2}\pmod{\Phi_n(q)},
\end{align*}
for $1\leq j\leq n$.
Similarly, there holds
\begin{align*}
&\hskip -2mm \qbinom{-1-a-sn}{k}{q^d}\\
&\equiv
\qbinom{-1-a}{k}{q^d}-\sum_{j=1}^{k}(-1)^{j}q^{-dj(k-j)-d\binom{j}2}\frac{[-sn]_{q^d}}{[j]_{q^d}}\qbinom{-1-a}{k-j}{q^d}\pmod{\Phi_n(q)^2}.
\end{align*}
Therefore,
\begin{align}
&\hskip -2mm \sum_{k=0}^{n-1}q^{dk^2}\qbinom{a+sn}{k}{q^d}\qbinom{-1-a-sn}{k}{q^d}-\sum_{k=0}^{n-1}q^{dk^2}\qbinom{a}{k}{q^d}\qbinom{-1-a}{k}{q^d}
\nonumber \\[5pt]
&\equiv -\sum_{k=1}^{n-1}q^{dk^2}\qbinom{a}{k}{q^d}\sum_{j=1}^{k}(-1)^{j}q^{-dj(k-j)-d\binom{j}2}\frac{[sn]_{q^d}}{[j]_{q^d}}\qbinom{-1-a}{k-j}{q^d}
\nonumber \\[5pt]
&\quad{}-\sum_{k=1}^{n-1}q^{dk^2}\qbinom{-1-a}{k}{q^d}
\sum_{j=1}^{k}(-1)^{j}q^{-dj(k-j)-d\binom{j}2}\frac{[sn]_{q^d}}{[j]_{q^d}}\qbinom{a}{k-j}{q^d}\pmod{\Phi_n(q)^2}. \label{eq:final-1}
\end{align}

By the $q$-Chu-Vandemonde identity \eqref{eq:qChu-2}, we have
\begin{align}
&\sum_{k=1}^{n-1}q^{dk^2}\qbinom{a}{k}q
\sum_{j=1}^{k}\frac{(-1)^{j}q^{-dj(k-j)-d\binom{j}2}}{[j]_{q^d}}\qbinom{-1-a}{k-j}{q^d}
\nonumber \\[5pt]
&=\sum_{j=1}^{a}\frac{(-1)^{j}q^{dj^2-d\binom{j}2}}{[j]_{q^d}}\sum_{k=j}^{a}q^{dk(k-j)}\qbinom{a}{a-k}{q^d}
\qbinom{-1-a}{k-j}{q^d}  \nonumber  \\[5pt]
&=\sum_{j=1}^{a}\frac{(-1)^{j}q^{dj^2-d\binom{j}2}}{[j]_{q^d}}\qbinom{-1}{a-j}{q^d}  \nonumber  \\[5pt]
%&\equiv(-1)^a\sum_{j=1}^{a}\frac{q^{dj^2-d\binom{j}2-d\binom{a-j+1}{2}}}{[j]_{q^d}}\\
&=(-1)^a\sum_{j=1}^{a}\frac{q^{-\frac{d(a+1)(a-2j)}{2}}}{[j]_{q^d}}. \label{eq:final-2}
\end{align}
Similarly, we have
\begin{align}
&\hskip -2mm \sum_{k=1}^{n-1}q^{dk^2}\qbinom{-1-a}{k}{q^d}
\sum_{j=1}^{k}\frac{(-1)^{j}q^{-dj(k-j)-d\binom{j}2}}{[j]_{q^d}}\qbinom{a}{k-j}{q^d}
\nonumber \\[5pt]
&\equiv \sum_{k=1}^{n-1}q^{dk^2}\qbinom{n-1-a}{k}{q^d}
\sum_{j=1}^{k}\frac{(-1)^{j}q^{-dj(k-j)-d\binom{j}2}}{[j]_{q^d}}\qbinom{1-(n-1-a)}{k-j}{q^d}
\nonumber \\[5pt]
&\equiv (-1)^{n-1-a}\sum_{j=1}^{n-1-a}\frac{q^{-\frac{d(n-a)(n-1-a-2j)}{2}}}{[j]_{q^d}} \nonumber \\[5pt]
&\equiv
(-1)^{n-1-a}\sum_{j=a+1}^{n-1}\frac{q^{-\frac{d(n-a)(2j-n-1-a)}{2}}}{[n-j]_{q^d}}\pmod{\Phi_n(q)}. \label{eq:final-3}
\end{align}
If $n$ is even, then
$$q^{\frac{n}2}=-1+\frac{1-q^n}{1-q^{\frac{n}{2}}}\equiv-1\pmod{\Phi_n(q)},$$
and (since $d$ is odd in this case)
$$
q^{-\frac{d(n-a)(2j-n-1-a)}{2}}\equiv q^{\frac{da(2j-1-a)-dn(2j-1)}{2}}\equiv-q^{\frac{da(2j-1-a)}{2}}\pmod{\Phi_n(q)};
$$
while if $n$ is odd, then
$$
q^{-\frac{d(n-a)(2j-n-1-a)}{2}}= q^{\frac{da(2j-1-a)-dn(2j-1-n)}{2}}
\equiv q^{\frac{da(2j-1-a)}{2}}\pmod{\Phi_n(q)}.
$$
Hence, we always have
\begin{align}
(-1)^{n-1-a}\sum_{j=a+1}^{n-1}\frac{q^{-\frac{d(n-a)(2j-n-1-a)}{2}}}{[n-j]_{q^d}}\equiv
%-\sum_{j=a+1}^{n-1}\frac{q^{dj+\frac{da(2j-1-a)}{2}}}{[j]_{q^d}}=
(-1)^{a-1}\sum_{j=a+1}^{n-1}\frac{q^{-\frac{d(a+1)(a-2j)}{2}}}{[j]_{q^d}}\pmod{\Phi_n(q)}.  \label{eq:final-4}
\end{align}

Noticing that
\begin{align*}
\sum_{j=1}^{n-1}\frac{1}{[j]_{q^d}}
&=\frac12\sum_{j=1}^{n-1}\bigg(\frac{1}{[j]_{q^d}}+\frac{1}{[n-j]_{q^d}}\bigg)\\
&=\frac12\sum_{j=1}^{n-1}\bigg(\frac{1}{[j]_{q^d}}+\frac{q^{jd}}{[n]_{q^d}-[j]_{q^d}}\bigg)\\
&\equiv
\frac12\sum_{j=1}^{n-1}\frac{1-q^{jd}}{[j]_{q^d}}=\frac{n-1}{2}(1-q^d)\pmod{\Phi_n(q)},
\end{align*}
we have
\begin{align}
\sum_{j=1}^{n-1}\frac{q^{d(a+1)j}}{[j]_{q^d}}
&=\sum_{j=1}^{n-1}\frac{1}{[j]_{q^d}}-(1-q^d)\sum_{j=1}^{n-1}\sum_{k=0}^{a}q^{dkj} \nonumber \\
&=\sum_{j=1}^{n-1}\frac{1}{[j]_{q^d}}-(n-1)(1-q^d)-(1-q^d)\sum_{k=1}^{a}\frac{q^{dk}-q^{dkn}}{1-q^{dk}} \nonumber\\
&\equiv \frac{n-1}{2}(1-q^d)-(n-1)(1-q^d)+a(1-q^d)\pmod{\Phi_n(q)}.  \label{eq:final-5}
\end{align}
It follows from \eqref{eq:final-1}--\eqref{eq:final-5} that
\begin{align*}
&\hskip -2mm
\sum_{k=0}^{n-1}q^{dk^2}\qbinom{a+sn}{k}{q^d}\qbinom{-1-a-sn}{k}{q^d}\\
&\equiv\sum_{k=0}^{n-1}q^{dk^2}\qbinom{a}{k}{q^d}\qbinom{-1-a}{k}{q^d}+\frac{2a+1-n}{2}(-1)^{a}q^{-d\binom{a+1}2}(1-q^{sdn}) \\
&=\qbinom{-1}{a}{q^d}+\frac{2a+1-n}{2}(-1)^{a}q^{-d\binom{a+1}2}(1-q^{sdn}) \pmod{\Phi_n(q)^2}
\end{align*}
by the $q$-Chu-Vandemonde identity \eqref{eq:qChu-2}.
Noticing that $\qbinom{-1}{a}{q^d}=(-1)^aq^{-d\binom{a+1}2}$ and
\begin{align*}
q^{(\alpha-a)d(\frac{n-1}2-a)}
&=(1-(1-q^{sdn}))^{\frac{n-1}{2}-a} \\
&\equiv 1+\frac{2a+1-n}{2}(1-q^{sdn})\pmod{\Phi_n(q)^2},
\end{align*}
we complete the proof.
\qed

\end{document}